\documentclass[12pt]{article}
\usepackage{amssymb,amsfonts,amsmath,graphicx}

\def\0{\leqno}
\def\({\left(}
\def\){\right)}
\def\<{\left<}
\def\>{\right>}
\def\ov{\overline}
\def\4{\subseteq }
\def\dd{\displaystyle}

\def\bit{\begin{itemize}}
\def\eit{\end{itemize}}
\def\barr{\begin{array}}
\def\earr{\end{array}}

\def\X#1#2{\stb{#1}{#2}{\mbox{\Huge$\times$}}}
\def\Z{{\rlap{$\kern2pt{\rm Z}$}{\rm Z}\,}}
\def\bld#1#2{{\buildrel{#1}\over{#2}}}
\def\st#1#2{{\mathrel{\mathop{#2}\limits_{#1}}{}\!}}
\def\stb#1#2#3{{\st{{#1}}{\bld{{#2}}{#3}}{}\!}}

\def\dd{\displaystyle}

\title{\bf A note on fundamental group lattices}
\author{Marius T\u arn\u auceanu}
\date{October 1, 2012}

\begin{document}

\maketitle

\begin{abstract}
The main goal of this note is to provide a new proof of a
classical result about projectivities between finite abelian
groups. It is based on the concept of fundamental group lattice,
studied in our previous papers \cite{8} and \cite{9}. A
generalization of this result is also given.
\end{abstract}



\section{Introduction}

The relation between the structure of a group and the structure of
its subgroup lattice constitutes an important domain of research
in group theory. One of the most interesting problems concerning
it is to study whether a group $G$ is determined by the subgroup
lattice of the $n$-th direct power $G^n$, $n\in\mathbb{N}^*$. In
other words, if the $n$-th direct powers of two groups have
isomorphic subgroup lattices, are these groups isomorphic? For
$n=1$ it is well-known that this problem has a negative answer
(see \cite{4}). The same thing can be also said for $n=2$,
excepting for some particular classes of groups, as simple groups
(see \cite{5}), finite abelian groups (see \cite{3}) or abelian
groups with the square root property (see \cite{2}). In the
general case (when $n\geq2$ is arbitrary) we recall Remark 1 of
\cite{2}, which states that an abelian group is determined by the
subgroup lattice of its $n$-th direct power if and only if it has
the $n$-th root property. This follows from some classical results
of Baer \cite{1}.

The starting point for our discussion is given by the papers
\cite{8} and \cite{9} (see also Section I.2.1 of \cite{7}), where
the concept of fundamental group lattice is introduced and
studied. It gives an arithmetic description of the subgroup
lattice of a finite abelian group and has many applications.
Fundamental group lattices were successfully used in \cite{8} to
solve the problem of existence and uniqueness of a finite abelian
group whose subgroup lattice is isomorphic to a fixed lattice and
in \cite{9} to count some types of subgroups of a finite abelian
group. In this paper they will be used to prove that the finite
abelian groups are determined by the subgroup lattices of their
direct $n$-powers, for any $n\geq2$. Notice that our proof is more
simple as the original one. A more general result will be also
inferred.

Most of our notation is standard and will usually not be repeated
here. Basic definitions and results on groups can be found in
\cite{6}. For subgroup lattice notions we refer the reader to
\cite{4} and \cite{7}.
\bigskip

In the following we recall the concept of fundamental group
lattice and two related theorems. Let $G$ be a finite abelian
group and $L(G)$ be the subgroup lattice of $G$. Then, by the
fundamental theorem of finitely ge\-ne\-ra\-ted abelian groups,
there exist (uniquely determined by $G$) the numbers
$k\in\mathbb{N}^*$ and
$d_1,d_2,...,d_k\in\mathbb{N}\setminus\{0,1\}$ satisfying
$d_1|\hspace{0,5mm}d_2|...|\hspace{0,5mm}d_k$ such that
$$G\cong\X{i=1}{k}\mathbb{Z}_{d_i}.\0(*)$$This decomposition of a $G$
into a direct product of cyclic groups together with the form of
subgroups of $\mathbb{Z}^k$ (see Lemma 2.1 of \cite{8}) lead us to
the following construction:

Let $k \geq 1$ be an integer. Then, for every
$(d_1,d_2,...,d_k)\in(\mathbb{N}{\setminus}\{0,1\})^k,$ we
consider the set $L_{(k;d_1,d_2,...,d_k)}$ consisting of all
matrices $A=(a_{ij})\in{\cal M}_k(\mathbb{Z})$ which satisfy the
conditions:

\bigskip
$\begin{tabular}{rl}
I.& $a_{ij}=0,$ for any $i>j,$\vspace*{1,2mm}\\
II.& $0\le a_{1j},a_{2j},...,a_{j-1j}<a_{jj},$ for any
 $j=\overline{1,k}$,\vspace*{1,2mm}\\
III.& 1)\ $a_{11}|\hspace{0,5mm}d_1,$\vspace*{1,2mm}\\
 &2)\ $a_{22}|\Bigl(d_2,d_1\,\displaystyle\frac{a_{12}}{a_{11}}\Bigr),$\end{tabular}$

$\begin{tabular}{rl}
 \ \ \ \ \ \ &3)\
 $a_{33}|\Bigl(d_3,d_2\,\displaystyle\frac{a_{23}}{a_{22}},d_1\,\displaystyle\frac{\left|\begin{array}{ll}a_{12}&a_{13}\\
 a_{22}&a_{23}\end{array}\right|}{a_{22}a_{11}}\Bigr),$\\
 &$\vdots$\\
 &k)\ $a_{kk}|\Bigl(d_k,d_{k-1}\,\displaystyle\frac{a_{k-1k}}{a_{k-1\,k-1}},d_{k-2}\,
 \displaystyle\frac{\left|\begin{array}{ll}a_{k-2\,k-1}&a_{k-2k}\\
 a_{k-1\,k-1}&a_{k-1k}\end{array}\right|}
 {a_{k-1\,k-1}a_{k-2\,k-2}},...,$\vspace*{4,5mm}\\
 &\hfill$d_1\,\displaystyle\frac{\left|\begin{array}{llcl}
 a_{12}&a_{13}&\cdots&a_{1k}\\
 a_{22}&a_{23}&\cdots&a_{2k}\\
 \vdots&\vdots&&\vdots\vspace*{-3mm}\\
 0&0&\cdots&a_{k-1\,k}\end{array}\right|}{a_{k-1\,k-1}a_{k-2\,k-2}...
 a_{11}}\Bigr),$
 \end{tabular}$
\medskip

\noindent where by $(x_1,x_2,...,x_m)$ we denote the greatest
common divisor of the numbers $x_1,x_2,...,x_m\in\mathbb{Z}$. On
the set $L_{(k;d_1,d_2,...,d_k)}$ we introduce the ordering
relation $"\le"$, defined as follows: for $A=(a_{ij})$,
$B=(b_{ij})\in L_{(k;d_1,d_2,...,d_k)}$, put $A\le B$ if and only
if we have
\begin{itemize}
 \item[$1)'$] $b_{11}|\hspace{0,5mm}a_{11},$\vspace*{-1,5mm}
 \item[$2)'$] $b_{22}|\Bigl(a_{22},\displaystyle\frac{\left|\begin{array}{ll}
 a_{11}&a_{12}\\
 b_{11}&b_{12}
 \end{array}\right|}{b_{11}}\Bigr),$\vspace*{-1,5mm}
 \item[$3)'$] $b_{33}|\Bigl(a_{33},\displaystyle\frac{\left|\begin{array}{ll}
 a_{22}&a_{23}\\
 b_{22}&b_{23}\\
 \end{array}\right|}{b_{22}},
 \displaystyle\frac{\left|\begin{array}{lll}
 a_{11}&a_{12}&a_{13}\\
 b_{11}&b_{12}&b_{13}\\
 0&b_{22}&b_{23}
 \end{array}\right|}{b_{22}b_{11}}\Bigr),$\vspace*{-1,5mm}
 \item[$\vdots$]\vspace*{-3,5mm}
 \item[${\rm k})'$] $b_{kk}|\Bigl(a_{kk},\displaystyle\frac{\left|\begin{array}{ll}
 a_{k-1\,k-1}&a_{k-1\,k}\\
 b_{k-1\,k-1}&b_{k-1\,k}
 \end{array}\right|}{b_{k-1\,k-1}},
 \displaystyle\frac{\left|\begin{array}{lll}
 a_{k-2\,k-2}&a_{k-2\,k-1}&a_{k-2\,k}\\
 b_{k-2\,k-2}&b_{k-2\,k-1}&b_{k-2\,k}\\
 0&b_{k-1\,k-1}&b_{k-1\,k}
 \end{array}\right|}{b_{k-1\,k-1}b_{k-2\,k-2}},...,$\bigskip\\
 \ \hspace*{6cm}$\displaystyle\frac{\left|\begin{array}{llcl}
 a_{11}&a_{12}&\cdots&a_{1k}\\
 b_{11}&b_{12}&\cdots&b_{1k}\\
 \vdots&\vdots&&\vdots\\
 0&0&\cdots&b_{k-1\,k}\end{array}\right|}{b_{k-1\,k-1}b_{k-2\,k-2}...b_{11}}\Bigr).$\vspace*{-1,5mm}
\end{itemize}

\noindent Then $L_{(k;d_1,d_2,...,d_k)}$ forms a complete modular
lattice with respect to $\le$, called a {\it fundamental group
lattice of degree $k$}. A powerful connection between this lattice
and $L(G)$ has been established in \cite{8}.

\bigskip\noindent{\bf Theorem A.} {\it If $G$ is a finite abelian group with the
decomposition $(*)$, then its subgroup lattice $L(G)$ is
isomorphic to the fundamental group lattice
$L_{(k;d_1,d_2,...,d_k)}$.}
\bigskip

In order to study when two fundamental group lattices are
isomorphic (that is, when two finite abelian groups are
lattice-isomorphic), the following notation is useful. For every
integer $n\geq2$, we denote by $\pi(n)$ the set consisting of all
primes dividing $n$. Let
$d_i,d'_{i'}\in\mathbb{N}\setminus\{0,1\}$, $i=\ov{1,k}$,
$i'=\ov{1,k'}$, such that
$d_1|\hspace{0,5mm}d_2|...|\hspace{0,5mm}d_k$ and
$d'_1|\hspace{0,5mm}d'_2|...|\hspace{0,5mm}d'_{k'}$. Then we shall
write
$$(d_1,d_2,...,d_k)\sim(d'_1,d'_2,...,d'_{k'})$$whenever the
next three conditions are satisfied:
\begin{itemize}
\item[a)] $k=k'$.
\item[b)] $d_i=d'_i,\ i=\ov{1,k-1}$.
\item[c)] The sets $\pi(d_k)\setminus\pi\left(\dd\prod_{i=1}^{k-1}d_i\right)$ and
$\pi(d'_k)\setminus\pi\left(\dd\prod_{i=1}^{k-1}d'_i\right)$ have
the same number of elements, say $r$. Moreover, for $r=0$ we have
$d_k=d'_k$ and for $r \geq 1$, by denoting
$\pi(d_k)\setminus\pi\left(\dd\prod_{i=1}^{k-1}d_i\right)=\{p_1,p_2,...,p_r\},$
$\pi(d'_k)\setminus\pi\left(\dd\prod_{i=1}^{k-1}d'_i\right)=\{q_1,q_2,...,q_r\},$
we have
$$\dd\frac{d_k}{d'_k}=\prod_{j=1}^r\left(\dd\frac{p_j}{q_j}\right)^{s_j},$$where
$s_j\in\mathbb{N}^*$, $j=\overline{1,r}$.
\end{itemize}

The following theorem of \cite{8} will play an essential role in
proving our main results.

\bigskip\noindent{\bf Theorem B.} {\it Two fundamental group lattices $L_{(k;d_1,d_2,...,d_k)}$ and
$L_{(k';d'_1,d'_2,...,d'_{k'})}$ are isomorphic if and only if
$(d_1,d_2,...,d_k)\sim(d'_1,d'_2,...,d'_{k'})$.}

\section{Main results}

As we already have mentioned, there exist large classes of
non-isomorphic finite abelian groups whose lattices of subgroups
are isomorphic. Simple examples of such groups are easily obtained
by using Theorem B:
\begin{itemize}
\item[1.] $G=\mathbb{Z}_6 \mbox{ and } H=\mathbb{Z}_{10} \mbox{ (cyclic groups)},$
\item[2.] $G=\mathbb{Z}_2\times\mathbb{Z}_6 \mbox{ and } H=\mathbb{Z}_2\times\mathbb{Z}_{10} \mbox{ (non-cyclic groups)}.$
\end{itemize}
Moreover, Theorem B allows us to find a subclass of finite abelian
groups which are determined by their lattices of subgroups (see
also Proposition 2.8 of \cite{8}).

\bigskip\noindent{\bf Theorem 1.} {\it Let $G$ and $H$ be two finite abelian groups
such that one of them possesses a decomposition of type $(*)$ with
$\pi(d_k)=\pi\left(\dd\prod_{i=1}^{k-1}d_i\right)$. Then $G\cong
H$ if and only if $L(G)\cong L(H)$.}
\bigskip

Next we shall focus on isomorphisms between the subgroup lattices
of the direct $n$-powers of two finite abelian groups, for
$n\geq2$. An alternative proof of the following well-known result
can be also inferred from Theorem B.

\bigskip\noindent{\bf Theorem 2.} {\it Let $G$ and $H$ be two finite abelian groups.
Then $G\cong H$ if and only if $L(G^n)\cong L(H^n)$ for some
integer $n\geq 2$.}

\bigskip\noindent{\bf Proof.} Let
$G\cong\X{i=1}{k}\mathbb{Z}_{d_i}$ and
$H\cong\X{i=1}{k'}\mathbb{Z}_{d'_i}$ be the corresponding
decompositions $(*)$ of $G$ and $H$, respectively, and assume that
$L(G^n)\cong L(H^n)$ for some integer $n\geq 2$. Then the
fundamental group lattices
$$L_{(k;\hspace{1mm}\underbrace{d_1,d_1,...,d_1}_{n\ {\rm
factors}},...,\underbrace{d_k,d_k,...,d_k}_{n\ {\rm
factors}})}\mbox{ and
}L_{(k';\hspace{1mm}\underbrace{d'_1,d'_1,...,d'_1}_{n\ {\rm
factors}},...,\underbrace{d'_{k'},d'_{k'},...,d'_{k'}}_{n\ {\rm
factors}})}$$are isomorphic. By Theorem B, one obtains
$$(\underbrace{d_1,d_1,...,d_1}_{n\ {\rm
factors}},...,\underbrace{d_k,d_k,...,d_k}_{n\ {\rm
factors}})\sim(\underbrace{d'_1,d'_1,...,d'_1}_{n\ {\rm
factors}},...,\underbrace{d'_{k'},d'_{k'},...,d'_{k'}}_{n\ {\rm
factors}})$$and therefore $k=k'$ and $d_i=d'_i$, for all
$i=\ov{1,k}$. These equalities show that $G\cong H$, which
completes the proof.
\hfill\rule{1,5mm}{1,5mm}
\bigskip

Clearly, two finite abelian groups $G$ and $H$ satisfying
$L(G^m)\cong L(H^n)$ for some (possibly different) integers
$m,n\geq2$ are not necessarily isomorphic. Nevertheless, a lot of
conditions of this type can lead to $G\cong H$, as shows the
following theorem.

\bigskip\noindent{\bf Theorem 3.} {\it Let $G$ and $H$ be two finite abelian groups.
Then $G\cong H$ if and only if there are the integers $r\geq1$ and
$m_1, m_2,..., m_r, n_1,n_2,..., n_r\geq2$ such that
$(m_1,m_2,...,m_r)=(n_1,n_2,...,n_r)$ and $L(G^{m_i})\cong
L(H^{n_i})$, for all $i=\ov{1,r}$.}

\bigskip\noindent{\bf Proof.} Suppose that $G$ and $H$ have the
decompositions in the proof of Theorem 2. For every $i=1,2,...,r$,
the lattice isomorphism $L(G^{m_i})\cong L(H^{n_i})$ implies that
$km_i=k'n_i$, in view of Theorem B. Set $d=(m_1,m_2,...,m_r)$.
Then $d=\dd\sum_{i=1}^r\alpha_im_i$ for some integers
$\alpha_1,\alpha_2,...,\alpha_r$, which leads to
$$kd=k\dd\sum_{i=1}^r\alpha_im_i=\dd\sum_{i=1}^r\alpha_ikm_i=\dd\sum_{i=1}^r\alpha_ik'n_i=k'\dd\sum_{i=1}^r\alpha_in_i.$$Since
$d\hspace{0,5mm}|\hspace{0,5mm}n_i$, for all $i=\ov{1,r}$, we
infer that $k'\hspace{0,5mm}|\hspace{0,5mm}k$. In a similar manner
one obtains $k\hspace{0,5mm}|\hspace{0,5mm}k'$, and thus $k=k'$.
Hence $m_i=n_i$ and the group isomorphism $G\cong H$ is obtained
from Theorem 2.
\hfill\rule{1,5mm}{1,5mm}
\bigskip

Finally, we indicate an open problem concerning the above results.

\bigskip\noindent{\bf Open problem.} In Theorem 3 replace the
condition $(m_1,m_2,...,m_r)=(n_1,n_2,...,n_r)$ with other
connections between the numbers $m_i$ and $n_i$, $i=1,2,...,r$,
such that the respective equivalence be also true.

\vspace*{5ex}\small

\hfill
\begin{minipage}[t]{5cm}
Marius T\u arn\u auceanu \\
Faculty of  Mathematics \\
``Al.I. Cuza'' University \\
Ia\c si, Romania \\
e-mail: {\tt tarnauc@uaic.ro}
\end{minipage}

\end{document}